\documentclass[a4paper]{article}
\usepackage{amsfonts}
\usepackage{amssymb}
\usepackage{amsthm}
\usepackage{amsmath}
\usepackage{latexsym}
\usepackage[T1]{fontenc}
\usepackage[latin1]{inputenc}
\usepackage{color}

\def \eop {\hbox{}\nobreak\hfill
\vrule width 2mm height 2mm depth 0mm
\par \goodbreak \smallskip}

\def \eop {\hbox{}\nobreak\hfill \vrule width 2.0mm height 1.8mm depth 0mm
\par \goodbreak \smallskip}
\numberwithin{equation}{section}

\newtheorem{definition}{Definition}[section]
\newtheorem{theorem}{Theorem}[section]
\newtheorem{proposition}{Proposition}[section]
\newtheorem{lemma}{Lemma}[section]

\def \eop {\hbox{}\nobreak\hfill
\vrule width 2mm height 2mm depth 0mm
\par \goodbreak \smallskip}

\input{tcilatex}
\begin{document}

\author{Badreddine Mansouri, Mostapha abd el ouahab Saouli
\\
{University Mohamed Khider, PO BOX 145, 07000 Biskra, Algeria}}

\title{Reflected solutions of Anticipated Backward Doubly SDEs driven by
Teugels Martingales.}
\date{}
\maketitle


\renewcommand{\thefootnote}{
\fnsymbol{footnote}} \footnotetext{{\scriptsize E-mail addresses:
mansouri.badreddine@gmail.com
(Badreddine\ Mansouri), saoulimoustapha@yahoo.fr (Mostapha abd el ouahab Saouli)}}



\noindent \textbf{Abstract.} We deal with reflected solutions of anticipated backward
doubly stochastic differential equations (RABDSDEs) driven by Teugels
martingales associated with Lévy process under a Lipschitz generator where
the coefficients of these BDSDEs depend on the future and present value of
the solution $\left( Y,Z\right) $. Also we study the existence of a solution
for anticipated BDSDEs.

\bigskip


\section{\textbf{Introduction}}

\qquad Backward stochastic differential equations (BSDEs in short) were
introduced by Bismut for the linear case [2] and by Pardoux and Peng in the
general case [7]. Precisely, according to [7], given a data $(\xi ,f)$
consisting of a square integrable random variable $\xi $ and a progressively
measurable process $f$, a so-called generator, they proved the existence and
uniqueness of a solution to these equations. Recently, a new type of BSDE,
called anticipated BSDE (ABSDE in short), which can be regarded as a new
duality type of stochastic differential delay equations, was introduced by
Peng and Yang [9], see also [13, 14]. After they introduced the theory of BSDEs, Pardoux
and Peng in [8] considered a new kind of BSDEs, that is a class of backward
doubly stochastic dfferential equations (BDSDEs in short) with two different
directions of stochastic integrals. They proved existence and uniqueness of
solutions for BDSDEs under Lipschitz conditions on the coefficients.

Recently, a new type of BDSDE called anticipated BDSDE (ABDSDE in short),
which can be regarded as a new duality type of stochastic differential delay
equations, was introduced by Xu [16], The ABDSDE is of the form%
\begin{equation*}
\begin{array}{l}
Y_{t}=\xi+\int_{t}^{T}f(s,\Lambda _{s},\Lambda _{s}^{\phi ,\psi
})ds+\int_{t}^{T}g(s,\Lambda _{s},\Lambda _{s}^{\phi ,\psi })d\overleftarrow{%
B}_{s}-\int_{t}^{T}Z_{s}dW_{s},\text{ }t\in \left[ 0,T\right] , \\
\\
\left( Y_{t},Z_{t}\right) =\left( \eta _{t},\vartheta _{t}\right) ,\qquad
\qquad \qquad \qquad \qquad \qquad \qquad \qquad \qquad \qquad \qquad t\in
\left[ T,T+K\right] ,%
\end{array}%
\end{equation*}%
where, $\Lambda _{s}=\left( Y_{s},Z_{s}\right) $ and $\Lambda _{s}^{\phi
,\psi }=\left( Y_{s+\phi \left( s\right) },Z_{s+\psi \left( s\right)
}\right) $, $\phi \left( \cdot \right) :[0,T]\rightarrow
\mathbb{R}
_{+}\setminus \left\{ 0\right\} $ and $\psi \left( \cdot \right)
:[0,T]\rightarrow
\mathbb{R}
_{+}\setminus \left\{ 0\right\} $ are continuous functions and $f$ is called
a generator.

In the paper of Nulart et al. [5], a martingale representation theorem associated to Lévy processes was
proved. Then it is natural to extend BDSDEs driven by Brownian motion to
BDSDEs driven by a Lévy process. In the work of Ren et al. [10], the authors proved the existence
and uniqueness of solutions of BDSDEs driven by Teugels martingales
associated with a Lévy process, under Lipschitz conditions on the generator $
f$. These results were important from a pure mathematical point of view as
well as from an application point of view in the world of finance.

The first work of Reflected BDSDEsis itroduced by Bahlali et al. [1], after Y. Ren [10] introduced a special class of reflected BDSDEs
(RBDSDEs, in short), which is a BDSDE but the solution is forced to stay
above a lower barrier.

Motivated by the above results and by the result introduced by Xiaoming Xu
[16], we establish firstly the existence and uniqueness of the solution of
the Anticipated reflected BDSDE driven by Teugles Martingales (RABDSDEs, in
short) in the proof we using the result of Y. Ren [11]. Let us point out
that our papier extends the results of Y. Ren [11], Xiaoming Xu [16] and
Gaofeng Zong [17]. The main idea of the proof is to the point fixe theorem.
And we establish others results for ABDSDEs, we prove the existence and
uniqueness of the solution.

The organization of the paper is as follows. In Section 2, we give some
preliminaires on the martingales $\left\{ H^{\left( i\right) },t\geq
0\right\} $ and we consider the spaces of processus also we define the Itô's
lamma. In Section 3, under certain assumptions, we obtain the existence and
uniqueness solution for the associated Anticipated reflected BDSDEs
(RABDSDEs) and ABDSDEs.

\section{Preliminaries}

\qquad Let $\left( \Omega ,\mathcal{F},P,B_{t},L_{t};0\leq t\leq T\right) $
be a complete Brownien-Lévy space in $%
\mathbb{R}
\times
\mathbb{R}
^{\ast }$ with Lévy mesure. For $T>0$, let $\left\{ B_{t},0\leq t\leq
T\right\} $ is a standard Brownian motion defined on $\left( \Omega ,%
\mathcal{F},P\right) $ with values in $\mathbb{R}$ and $\left\{ L_{t};0\leq
t\leq T\right\} $ is a $\mathbb{R-}$valued pure jump-Lévy process of the
form $L_{t}=bt+l_{t}$ independent of $\left\{ B_{t},0\leq t\leq T\right\} $.

Let \ $\mathcal{F}_{t}^{L}:=\sigma (L_{s};0\leq s\leq t)$ \ and \ $\mathcal{F%
}_{t,T}^{B}:=\sigma (B_{s}-B_{t};t\leq s\leq T),$ completed with $P$-null
sets. We put,\ $\mathcal{F}_{t}:=\mathcal{F}_{t}^{L}\vee \mathcal{F}%
_{t,T}^{B}$, for each $t\in \left[ 0,T\right] $, and $\mathcal{H}_{t}:=%
\mathcal{F}_{0,t}^{L}\vee \mathcal{F}_{t,T+K}^{B}$ for each $t\in \left[
0,T+K\right] $. It should be noted that $\left( \mathcal{H}_{t}\right) $ is
not an increasing family of sub $\sigma -$fields, and hence it is not a
filtration.

For each $t\in \lbrack 0,T+K]$, we define%
\begin{equation*}
\mathcal{G}_{t}:=\mathcal{F}_{t}^{L}\vee \mathcal{F}_{T+K}^{B},
\end{equation*}%
the collection $\left( \mathcal{G}_{t}\right) _{t\in \left[ 0,T+K\right] }$
is a filtration.

For any $d,$ $k\geq 1$, we consider the following spaces of processus:

\begin{itemize}
\item Let $\mathcal{M}_{\mathcal{H}}^{2}\left( \left[ 0,T\right] ;\mathbb{R}%
\right) $ denote the set of \ $1-$dimensional, $\mathcal{H}_{t}-$%
progressively measurable stochastic processes $\left\{ \varphi _{t};t\in %
\left[ 0,T\right] \right\} $, such that \ $\mathbb{E}\int_{0}^{T}\left\vert
\varphi _{t}\right\vert ^{2}dt<\infty $.

\item We denote by $\mathcal{S}_{\mathcal{H}}^{2}\left( \left[ 0,T\right] ;%
\mathbb{R}\right) $, the set of continuous and $\mathcal{H}_{t}-$%
progressively measurable stochastic processes $\left\{ \varphi _{t};t\in %
\left[ 0,T\right] \right\} $, which satisfy $\mathbb{E}(\sup_{0\leq t\leq
T}\left\vert \varphi _{t}\right\vert ^{2})<\infty $.

\item $l^{2}$ be the space of real valued sequences $\left( x_{n}\right)
_{n\geq 0}$ such that $\sum_{i=1}^{i=\infty }x_{i}^{2}<\infty ,$ and $%
\left\vert \left\vert x\right\vert \right\vert
_{l^{2}}^{2}=\sum_{i=1}^{i=\infty }x_{i}^{2}.$

\item $\mathcal{A}^{2}$ set of continuous, increasing, $\mathcal{H}_{t}-$%
measurable process $K:\left[ 0,T\right] \times \Omega \rightarrow \lbrack
0,+\infty ($ with $K_{0}=0,$ $\mathbb{E}\left( K_{T}\right) ^{2}<+\infty .$

\item $\mathcal{M}_{\mathcal{H}}^{2}\left( \left[ 0,T\right] ;l^{2}\right) $
and $\mathcal{S}_{\mathcal{H}}^{2}\left( \left[ 0,T\right] ;l^{2}\right) $;
are the corresponding spaces of $l^{2}$-valued processes equipped with the
norm $\left\vert \left\vert \varphi \right\vert \right\vert _{l^{2}}^{2}=%
\mathbb{E}\int_{0}^{T}\sum_{i=1}^{i=\infty }\left\vert \varphi _{t}^{\left(
i\right) }\right\vert ^{2}dt<\infty .$

\item $\mathbb{L}^{2}\left( \mathcal{H}_{T}\right) $ set of $\mathcal{H}_{T}$%
- measurable random variables $\xi :\Omega \rightarrow
\mathbb{R}
^{d}$ with $\mathbb{E}\left\vert \xi \right\vert ^{2}<+\infty .$

\item Notice that the space $\mathcal{B}_{\mathcal{H}}^{2}\left( \left[ 0,T%
\right] ,%
\mathbb{R}
\right) =\mathcal{S}_{\mathcal{H}}^{2}\left( \left[ 0,T\right] ;\mathbb{R}%
\right) \times \mathcal{M}_{\mathcal{H}}^{2}\left( \left[ 0,T\right]
;l^{2}\right) $ endowd with the norm%
\begin{equation*}
\left\vert \left\vert \left( Y,Z\right) \right\vert \right\vert _{\mathcal{B}%
_{\mathcal{H}}^{2}\left( \left[ 0,T\right] ,%
\mathbb{R}
\right) }=\left\vert \left\vert Y\right\vert \right\vert _{\mathcal{S}_{%
\mathcal{H}}^{2}\left( \left[ 0,T\right] ;\mathbb{R}\right) }+\left\vert
\left\vert Z\right\vert \right\vert _{\mathcal{M}_{\mathcal{H}}^{2}\left( %
\left[ 0,T\right] ;l^{2}\right) }.
\end{equation*}
\end{itemize}

We denote by $\left( H^{i}\right) _{i\geq 1}$ the Teugels martingale
associated with the lévy process $\left\{ L_{t},t\in \left[ 0,T\right]
\right\} $ with is given by $H^{\left( i\right) }=c_{i,i}Y_{t}^{\left(
i\right) }+c_{i,i-1}Y_{t}^{\left( i-1\right) }+...+c_{i,1}Y_{t}^{\left(
1\right) },$ where $Y_{t}^{\left( i\right) }=L_{t}^{\left( i\right) }-%
\mathbb{E}\left[ L_{t}^{\left( i\right) }\right] =L_{t}^{\left( i\right) }-t%
\mathbb{E}\left[ L_{1}^{\left( i\right) }\right] $ for all $i\geq 1$ and $%
L_{t}^{\left( i\right) }$ are power-jump processes. That is, $L_{t}^{\left(
1\right) }=L_{t}$ and $L_{t}^{\left( i\right) }=\sum_{0\leq s\leq t}\left(
\Delta L_{s}\right) ^{i}$ for $i\geq 2,$ and $\left[ H^{\left( i\right)
},H^{\left( j\right) }\right] ,i\neq j$ and $\left\{ \left[ H^{\left(
i\right) },H^{\left( i\right) }\right] _{t}-t,\text{ }t\geq 0\right\} $ are
uniformly integrable martingale with initial value $0,$ i.e.,%
\begin{equation*}
\langle H^{\left( i\right) },H^{\left( j\right) }\rangle _{t}=t\delta _{i,j},
\end{equation*}%
where it was shown in [6] that the coefficients $c_{i,k}$ correspond to the
orthonormalization of the polynomials $1,x,x^{2},...$with respect to the
measure $\mu \left( dx\right) =x^{2}v\left( dx\right) +\sigma ^{2}\delta
_{0}(dx).$ the resulting processes $H^{\left( i\right) }=\left\{ H^{\left(
i\right) },t\geq 0\right\} $ are called the orthonormalized ith-power-jump
processes.

The result depends on the following extension of the well-krown Itô's
formula. Its proof follows the same way as lemma 1.3 of [8]

\begin{lemma}
\label{Lemma 1} Let $\alpha \in \mathcal{S}_{\mathcal{H}}^{2}\left( \left[
0,T\right] ;\mathbb{R}\right) ,$ $\beta ,$ $\gamma $ and $\sigma \in
\mathcal{M}_{\mathcal{H}}^{2}\left( \left[ 0,T\right] ;l^{2}\right) $ such
that%
\begin{equation*}
\alpha _{t}=\alpha _{0}+\int_{0}^{t}\beta _{s}ds+\int_{0}^{t}\gamma
_{s}dB_{s}+\sum_{i=1}^{\infty }\int_{0}^{t}\sigma _{s}^{\left( i\right)
}dH_{s}^{\left( i\right) },
\end{equation*}%
then%
\begin{eqnarray*}
\left\vert \alpha _{t}\right\vert ^{2} &=&\left\vert \alpha _{0}\right\vert
^{2}+2\int_{0}^{t}\alpha _{s}\beta _{s}ds+2\int_{0}^{t}\alpha _{s}\gamma
_{s}dB_{s}+2\sum_{i=1}^{\infty }\int_{0}^{t}\alpha _{s}\sigma _{s}^{\left(
i\right) }dH_{s}^{\left( i\right) } \\
&&-\int_{0}^{t}\left\vert \gamma _{s}\right\vert ^{2}ds+\sum_{i=1}^{\infty
}\sum_{j=1}^{\infty }\int_{0}^{t}\sigma _{s}^{\left( i\right) }\sigma
_{s}^{\left( j\right) }d\left[ H^{\left( i\right) },H^{\left( j\right) }%
\right] _{s},
\end{eqnarray*}%
note that $\langle H^{\left( i\right) },H^{\left( j\right) }\rangle
_{t}=\delta _{ij}t,$ we have%
\begin{equation*}
\mathbb{E}\left\vert \alpha _{t}\right\vert ^{2}\leq \mathbb{E}\left\vert
\alpha _{0}\right\vert ^{2}+2\mathbb{E}\int_{0}^{t}\alpha _{s}\beta _{s}ds-%
\mathbb{E}\int_{0}^{t}\left\vert \gamma _{s}\right\vert ^{2}ds+\mathbb{E}%
\sum_{i=1}^{\infty }\int_{0}^{t}\left( \sigma _{s}^{\left( i\right) }\right)
^{2}ds.
\end{equation*}%
%
%
%
%
%
%
%
%
%
%
%
%
%
%
%
%
\end{lemma}

\section{Main result}

\subsection{Anticipated BDSDE with Lower barrier.}

\qquad In this subsection, we consider the following $1-$dimensional
anticipated reflected backward doubly stochastic differential equation%
\begin{equation}
\left\{
\begin{array}{l}
Y_{t}=\xi+\int_{t}^{T}f(s,\Lambda _{s},\Lambda _{s}^{\phi ,\psi
})ds+\int_{t}^{T}g(s,\Lambda _{s},\Lambda _{s}^{\phi ,\psi })d\overleftarrow{%
B}_{s} \\
\\
+\int_{t}^{T}dK_{s}-\sum_{i=1}^{\infty }\int_{t}^{T}Z_{s}^{\left( i\right)
}dH_{s}^{\left( i\right) },\qquad \qquad \qquad \quad t\in \left[ 0,T\right]
, \\
\\
\left( Y_{t},Z_{t}\right) =\left( \eta _{t},\vartheta _{t}\right) ,\qquad
\qquad \qquad \qquad \qquad \qquad \text{\ \ }t\in \left[ T,T+K\right] ,%
\end{array}%
\right.
\end{equation}%
where $f$ called the generator, $\Lambda _{s}=\left( Y_{s-},Z_{s}\right) $
and $\Lambda _{s}^{\phi ,\psi }=\left( Y_{s+\phi \left( s\right)
-},Z_{s+\psi \left( s\right) }\right) .$

\vskip0.15cm Let $\phi :\left[ 0,T\right] \rightarrow
\mathbb{R}
_{+}^{\ast }$, and $\psi :\left[ 0,T\right] \rightarrow
\mathbb{R}
_{+}^{\ast }$ are continuous functions satisfying:\newline
\textbf{(A)} There exists a constant $K\geq 0$ such that for all $t\in \left[
0,T\right] ,$%
\begin{equation*}
t+\phi \left( t\right) \leq T+K,\qquad \qquad t+\psi \left( t\right) \leq
T+K.
\end{equation*}%
\textbf{(B)} There exists a constant $M\geq 0$ such that for each $t\in %
\left[ 0,T\right] $ and for all nonnegative integrable functions $h\left(
\cdot \right) $,%
\begin{equation*}
\left\{
\begin{array}{c}
\int_{t}^{T}h(s+\phi \left( s\right) )ds\leq M\int_{t}^{T+K}h(s)ds, \\
\\
\int_{t}^{T}h(s+\psi \left( s\right) )ds\leq M\int_{t}^{T+K}h(s)ds.%
\end{array}%
\right.
\end{equation*}

\begin{definition}
A solution of equation $(3.1)$ is a triple $(Y,Z,K)$ which belongs to the
space $\mathcal{B}_{\mathcal{H}}^{2}([0,\,T+K],\mathbb{R})\times \mathcal{A}%
^{2}$ and satisfies $(3.1)$ sach that:%
\begin{equation*}
\left\{
\begin{array}{c}
S_{t}\leq Y_{t},\text{ }0\leq t\leq T+K, \\
\int_{0}^{T}\left( Y_{s-}-S_{s-}\right) dK_{s}=0.%
\end{array}%
\right.
\end{equation*}
\end{definition}

In this subsection we study the ABDSDEs with reflection under Lipschitz
continuous generator. We consider the following assumptions \textbf{(H1)}:

\vskip0.1cm\noindent \textbf{(H1.1) }\noindent (i) There exist a constant $%
c>0$ such that for any $\left( r,\acute{r}\right) \in \left[ 0,T+K\right]
^{2}$, $\left( t,\omega ,y,z,\pi ,\zeta \right) ,$ $\left( t,\omega ,y,z,%
\acute{\pi},\acute{\zeta}\right) \in \left[ 0,T\right] \times \Omega \times
\mathbb{R}
\times l^{2}\times \mathcal{B}_{\mathcal{H}}^{2}\left( \left[ 0,T+K\right] ,%
\mathbb{R}
\right) ,$

\begin{eqnarray*}
&&\left\vert f(t,\omega ,y,z,\pi \left( r\right) ,\zeta \left( \acute{r}%
\right) )-f(t,\omega ,\acute{y},\acute{z},\acute{\pi}\left( r\right) ,\acute{%
\zeta}\left( \acute{r}\right) )\right\vert ^{2} \\
&\leq &c\left( \left\vert y-\acute{y}\right\vert ^{2}+\left\vert \left\vert
z-\acute{z}\right\vert \right\vert _{l^{2}}^{2}+\mathbb{E}^{\mathcal{F}_{t}}%
\left[ \left\vert \pi \left( r\right) -\acute{\pi}\left( r\right)
\right\vert ^{2}+\left\vert \left\vert \zeta \left( \acute{r}\right) -\acute{%
\zeta}\left( \acute{r}\right) \right\vert \right\vert _{l^{2}}^{2}\right]
\right) .
\end{eqnarray*}%
\vskip0.1cm\noindent (ii) There exists a constant $c>0$, $0<\alpha _{1}<%
\frac{1}{2}$ and $0<\alpha _{2}<\frac{1}{M}$ satisfying $0<\alpha
_{1}+\alpha _{2}M<\frac{1}{2}$, such that%
\begin{eqnarray*}
&&\left\vert g(t,\omega ,y,z,\pi \left( r\right) ,\zeta \left( r\right)
)-g(t,\omega ,\acute{y},\acute{z},\acute{\pi}\left( r\right) ,\acute{\zeta}%
\left( \acute{r}\right) )\right\vert ^{2} \\
&\leq &c\left( \left\vert y-\acute{y}\right\vert ^{2}+\mathbb{E}^{\mathcal{F}%
_{t}}\left\vert \pi \left( r\right) -\acute{\pi}\left( r\right) \right\vert
^{2}\right) +\alpha _{1}\left\vert \left\vert z-\acute{z}\right\vert
\right\vert _{l^{2}}^{2}+\alpha _{2}\mathbb{E}^{\mathcal{F}_{t}}\left\vert
\left\vert \zeta \left( \acute{r}\right) -\acute{\zeta}\left( \acute{r}%
\right) \right\vert \right\vert _{l^{2}}^{2}.
\end{eqnarray*}%
\vskip0.1cm\noindent \textbf{(H1.2) }For any\textbf{\ }$\left( t,\omega
,y,z,\pi ,\zeta \right) ,$%
\begin{eqnarray*}
\mathbb{E}\int_{0}^{T}\left\vert f(s,\omega ,0,0,0,0)\right\vert ds
&<&\infty , \\
\mathbb{E}\int_{0}^{T}\left\vert g(s,\omega ,y,z,\pi ,\zeta )\right\vert ds
&<&\infty .
\end{eqnarray*}%
\vskip0.1cm\noindent \textbf{(H1.3) }The terminal value\textbf{\ }$\xi $ be
a given random variable in $\mathbb{L}^{2}$.

Also we consider the following assumptions \textbf{(H2)}:

\vskip0.1cm\noindent\ \textbf{(H2.1) }$\left( S_{t}\right) _{t\geq 0},$ is a
continuous progressively measurable real valued process satisfying%
\begin{equation*}
\begin{tabular}{lll}
$\mathbb{E}\left( \sup_{0\leq t\leq T+K}\left( S_{t}^{+}\right) ^{2}\right)
<+\infty ,$ & where & $S_{t}^{+}:=\max \left( S_{t},0\right) .$%
\end{tabular}%
\end{equation*}%
\vskip0.1cm\noindent \textbf{(H2.2) }For any\textbf{\ }$t\in \left[ T,T+K%
\right] ,$ $S_{t}\leq \eta _{t}$, $\mathbb{P}$-almost surely.

\vskip0.1cm\noindent \textbf{(H2.3) }$\left( \eta _{t},\vartheta _{t}\right)
\in \mathcal{B}_{\mathcal{H}}^{2}\left( \left[ T,T+K\right] ,%
\mathbb{R}
\right) .$

\vskip0.1cm\noindent \textbf{(H2.4) }$( K_{t})_{t\in \left[
0,T\right]} $ is a continuous, increasing process with $K_{0}=0$ and $\mathbb{%
E}\left( K_{T}\right) ^{2}<+\infty .$

\subsubsection{\textbf{Existence and uniqueness of solutions.}}

\begin{theorem}
\label{theorem 3} Let $f,$ $g$ satisfies the hypothesis $\left( H1\right) $,
$\left( H2\right) $ and \textbf{(A), (B)} are hold. Then the Anticipated
RBDSDEs $\left( 3.1\right) $ has a unique solution $\left(
Y_{t},Z_{t},K_{t}\right) _{t\in \left[ 0,T+K\right] }.$
\end{theorem}

\subsubsection{Proof of the existence and uniqueness result.}

\noindent \textbf{Proof :} \textbf{Uniqueness.} Let $\left(
Y^{j},Z^{j},K^{j}\right) \in \mathcal{B}_{\mathcal{H}}^{2}([0,T+K],\mathbb{R}%
)\times \mathcal{A}^{2}$ for $j=1,2$ be any two solutions. Define $\Delta
Y_{t}=Y_{t}^{1}-Y_{t}^{2},$ $\Delta Z_{t}^{\left( i\right) }=Z_{t}^{1,\left(
i\right) }-Z_{t}^{2,\left( i\right) }$ and $\Delta K_{s}=K_{t}^{1}-K_{t}^{2}$
and for a function:%
\begin{equation*}
\left\{
\begin{array}{c}
\Delta f\left( s\right) =f(s,\Lambda _{s}^{1},\Lambda _{s}^{1,\phi ,\psi
})-f(s,\Lambda _{s}^{2},\Lambda _{s}^{2,\phi ,\psi }), \\
\\
\Delta g\left( s\right) =g(s,\Lambda _{s}^{1},\Lambda _{s}^{1,\phi ,\psi
})-g(s,\Lambda _{s}^{2},\Lambda _{s}^{2,\phi ,\psi }).%
\end{array}%
\right.
\end{equation*}%
We consider the following equation%
\begin{equation*}
\left\{
\begin{array}{l}
\Delta Y_{t}=\xi _{T}+\int_{t}^{T}\Delta f(s)ds+\int_{t}^{T}\Delta g(s)d%
\overleftarrow{B}_{s}+\int_{t}^{T}d\left( \Delta K_{s}\right)
-\sum_{i=1}^{\infty }\int_{t}^{T}\Delta Z_{s}^{\left( i\right)
}dH_{s}^{\left( i\right) },\ t\in \left[ 0,T\right] , \\
\\
\left( \Delta Y_{t},\Delta Z_{t}\right) =\left( 0,0\right) ,\qquad \qquad
\qquad \qquad \qquad \qquad \qquad \qquad \qquad \qquad \qquad \qquad \qquad
t\in \left[ T,T+K\right] ,%
\end{array}%
\right.
\end{equation*}%
where $\Lambda _{s}^{j}=\left( Y_{s-}^{j},Z_{s}^{j}\right) $ and $\Lambda
_{s}^{j,\phi ,\psi }=\left( Y_{s+\phi \left( s\right) -}^{j},Z_{s+\psi
\left( s\right) }^{j}\right) $ for $j=1,2.$

It follows from Itô's formula that%
\begin{eqnarray*}
\mathbb{E}\left( e^{\beta t}|\Delta Y_{t}|^{2}\right) &+&\beta \mathbb{E}%
\int_{t}^{T}e^{\beta s}\left\vert \Delta Y_{s}\right\vert
^{2}ds+\sum_{i=0}^{\infty }\mathbb{E}\int_{t}^{T}e^{\beta s}|\Delta
Z_{s}^{\left( i\right) }|^{2}ds \\
&=&2\mathbb{E}\int_{t}^{T}e^{\beta s}\Delta Y_{s}\Delta f\left( s\right) ds+%
\mathbb{E}\int_{t}^{T}e^{\beta s}|\Delta g\left( s\right) |^{2}ds+2\mathbb{E}%
\int_{t}^{T}e^{\beta s}\Delta Y_{s-}d\left( \Delta K_{s}\right) .
\end{eqnarray*}%
Since$\int_{t}^{T}e^{\beta s}\Delta Y_{s-}d\left( \Delta K_{s}\right) \leq
0, $ we have%
\begin{eqnarray*}
\mathbb{E}\left( e^{\beta t}|\Delta Y_{t}|^{2}\right) &+&\beta \mathbb{E}%
\int_{t}^{T}e^{\beta s}\left\vert \Delta Y_{s}\right\vert
^{2}ds+\sum_{i=0}^{\infty }\mathbb{E}\int_{t}^{T}e^{\beta s}|\Delta
Z_{s}^{\left( i\right) }|^{2}ds \\
&\leq &2\mathbb{E}\int_{t}^{T}e^{\beta s}\Delta Y_{s}\Delta f\left( s\right)
ds+\mathbb{E}\int_{t}^{T}e^{\beta s}|\Delta g\left( s\right) |^{2}ds.
\end{eqnarray*}%
Using Young's inequality $2ab\leq \epsilon _{1}a^{2}+\frac{b^{2}}{\epsilon
_{1}}$ and hypothesis $\left( H.1\right) ,$ we have%
\begin{equation*}
2\mathbb{E}\int_{t}^{T}e^{\beta s}\Delta Y_{s}\Delta f(s)ds\leq \epsilon _{1}%
\mathbb{E}\int_{t}^{T+K}e^{\beta s}\left\vert \Delta Y_{s}\right\vert
^{2}ds+\left( \frac{c}{\epsilon _{1}}+cM\right) \mathbb{E}%
\int_{t}^{T+K}e^{\beta s}\left( \left\vert \Delta Y_{s-}\right\vert
^{2}+\left\vert \left\vert \Delta Z_{s}\right\vert \right\vert
_{l^{2}}^{2}\right) ds,
\end{equation*}%
and also%
\begin{equation*}
\mathbb{E}\int_{t}^{T}e^{\beta s}\left\vert \Delta g(s)\right\vert
^{2}ds\leq \left( c+cM\right) \mathbb{E}\int_{t}^{T+K}e^{\beta s}\left\vert
\Delta Y_{s-}\right\vert ^{2}ds+\left( \alpha _{1}+\alpha _{2}M\right)
\mathbb{E}\int_{t}^{T+K}e^{\beta s}\left\vert \left\vert \Delta
Z_{s}\right\vert \right\vert _{l^{2}}^{2}ds.
\end{equation*}%
Then, we have the following inequality%
\begin{eqnarray*}
\mathbb{E}\left( e^{\beta t}|\Delta Y_{t}|^{2}\right) &+&\beta \mathbb{E}%
\int_{t}^{T}e^{\beta s}\left\vert \Delta Y_{s}\right\vert ^{2}ds+\mathbb{E}%
\int_{t}^{T}e^{\beta s}\left\vert |\Delta Z_{s}|\right\vert _{l^{2}}^{2}ds \\
&\leq &\epsilon _{1}\mathbb{E}\int_{t}^{T+K}e^{\beta s}\left\vert \Delta
Y_{s}\right\vert ^{2}ds+\left( \frac{c}{\epsilon _{1}}+2cM+c\right) \mathbb{E%
}\int_{t}^{T+K}e^{\beta s}\left\vert \Delta Y_{s-}\right\vert ^{2}ds \\
&&+\left( \alpha _{1}+\alpha _{2}M+\frac{c}{\epsilon _{1}}+cM\right) \mathbb{%
E}\int_{t}^{T+K}e^{\beta s}\left\vert \left\vert \Delta Z_{s}\right\vert
\right\vert _{l^{2}}^{2}ds,
\end{eqnarray*}%
chossing $\epsilon _{1}>0$ such that, $\left( \alpha _{1}+\alpha _{2}M+\frac{%
c}{\epsilon _{1}}+cM\right) <1$ and $\beta -\epsilon _{1}>0,$ we get%
\begin{equation*}
\mathbb{E}\left( e^{\beta t}|\Delta Y_{t}|^{2}\right) \leq C\mathbb{E}%
\int_{t}^{T+K}e^{\beta s}\left\vert \Delta Y_{s-}\right\vert ^{2}ds,
\end{equation*}%
where $C=\frac{c}{\epsilon _{1}}+2cM+c.$ The uniqueness of solution follows
from Gronwall's lemma.

\noindent \textbf{Existence.} Before we start proving equation $(3.1)$ has a unique solution with $f,$ $g$
independent on the value and the futur value of $\left( Y,Z\right) $, i.e.,
P-a.s.,$f\left( t,\omega ,y,z,\pi ,\zeta \right) =f\left( t,\omega \right) $
and $g\left( t,\omega ,y,z,\pi ,\zeta \right) =g\left( t,\omega \right) ,$
for any $\left( t,y,z,\pi ,\zeta \right) .$
Then by Y. Ren [11] and the provious proof, we deduce
that the l'equation $(3.1)$ where $f$ and $g$ independent on the value and
the futur value of $\left( Y,Z\right) $ has a unique solution.

Now, we shall prove the existence in the general case. For all $\left( r,%
\acute{r}\right) \in \left[ t,T+K\right] ^{2},$ $\forall t\in \left[ 0,T+K%
\right] $%
\begin{equation*}
f\left( t,\omega ,y,z,\pi \left( r\right) ,\zeta \left( \acute{r}\right)
\right) :\left[ 0,T\right] \times \Omega \times
\mathbb{R}
\times l^{2}\times \mathcal{B}_{\mathcal{H}}^{2}([0,T+K],\mathbb{R}%
)\rightarrow \mathcal{M}_{\mathcal{H}}^{2}([0,T+K],\mathbb{R}).
\end{equation*}%
Let $\mathcal{S}_{\mathcal{H}}^{2}\left( \left[ 0,T+K\right] ;\mathbb{R}%
^{d}\right) \times \mathcal{M}_{\mathcal{H}}^{2}\left( \left[ 0,T+K\right]
;l^{2}\right) $ endowed with the norm%
\begin{equation*}
\left\vert \left\vert \left( Y,Z\right) \right\vert \right\vert _{\beta
}=\left( \mathbb{E}\left[ \int_{0}^{T+K}e^{\beta s}\left( \left\vert
Y_{s-}\right\vert ^{2}ds+\sum_{i=1}^{i=\infty }\left\vert Z_{s}^{\left(
i\right) }\right\vert ^{2}\right) ds\right] \right) ^{\frac{1}{2}}.
\end{equation*}%
Therefore, for given $U\in \mathcal{S}_{\mathcal{H}}^{2}([0,T+K],\mathbb{R})$%
, $V\in \mathcal{M}_{\mathcal{H}}^{2}\left( \left[ 0,T+K\right]
;l^{2}\right) $, there exists a unique solution $\left( Y,Z,K\right) $ for
the following ABDSDEs with reflection%
\begin{equation}
\left\{
\begin{array}{l}
Y_{t}=\xi _{T}+\int_{t}^{T}f(s,\theta _{s},\theta _{s}^{\phi ,\psi
})ds+\int_{t}^{T}g(s,\theta _{s},\theta _{s}^{\phi ,\psi })d\overleftarrow{B}%
_{s}+\int_{t}^{T}dK_{s}-\sum_{i=1}^{\infty }\int_{t}^{T}Z_{s}^{\left(
i\right) }dH_{s}^{\left( i\right) },\ t\in \left[ 0,T\right] , \\
\\
\left( Y_{t},Z_{t}\right) =\left( \eta _{t},\vartheta _{t}\right) ,\qquad
\qquad \qquad \qquad \qquad \qquad \qquad \qquad \qquad \qquad \qquad \qquad
\quad \qquad \qquad t\in \left[ T,T+K\right] ,%
\end{array}%
\right.
\end{equation}%
We can construct the mapping $\Phi $ is well defined, let $\left(
Y_{t},Z_{t}\right) $ and $\left( \tilde{Y}_{t},\tilde{Z}_{t}\right) $ be two
solution of system $(3.2)$ such that $\left( Y_{t},Z_{t}\right) =\Phi \left(
U_{t-},V_{t}\right) $ and $\left( \tilde{Y}_{t},\tilde{Z}_{t}\right) =\Phi
\left( \tilde{U}_{t-},\tilde{V}_{t}\right) $.

For $\beta \in
\mathbb{R}
$. The couple $\left( \Delta Y_{t},\Delta Z_{t}\right) $ solve the ABDSDEs
with reflection%
\begin{equation*}
\left\{
\begin{array}{l}
\Delta Y_{t}=\int_{t}^{T}\Delta f(s)ds+\int_{t}^{T}\Delta g(s)d%
\overleftarrow{B}_{s}+\int_{t}^{T}d\left( \Delta K_{s}\right)
-\sum_{i=1}^{i=\infty }\int_{t}^{T}\Delta Z_{s}^{\left( i\right)
}dH_{s}^{\left( i\right) },\text{ }t\in \left[ 0,T\right] , \\
\\
\left( \Delta Y_{t},\Delta Z_{t}\right) =\left( 0,0\right) ,\qquad \qquad
\qquad \qquad \qquad \qquad \qquad \qquad \qquad \qquad \qquad \qquad t\in
\left[ T,T+K\right] .%
\end{array}%
\right.
\end{equation*}%
where for a function $h\in \left\{ f,g\right\} $, $\Delta h(s)=h(s,\theta
_{s},\theta _{s}^{\phi ,\psi })-h(s,\tilde{\theta}_{s},\tilde{\theta}%
_{s}^{\phi ,\psi }),$ $\theta _{s}=\left( U_{s-},V_{s}\right) ,$ $\theta
_{s}^{\phi ,\psi }=\left( U_{s+\phi \left( s\right) -},V_{s+\psi \left(
s\right) }\right) ,$ $\tilde{\theta}_{s}=\left( \tilde{U}_{s-},\tilde{V}%
_{s}\right) $, $\tilde{\theta}_{s}^{\phi ,\psi }=\left( \tilde{U}_{s+\phi
\left( s\right) -},\tilde{V}_{s+\psi \left( s\right) }\right) $ and $\Delta
\Psi _{s}=\Psi _{s}-\tilde{\Psi}_{s}$.

Now applying Itô's formula for $e^{\beta t}\left\vert \Delta
Y_{t}\right\vert ^{2},$ we get
\begin{eqnarray*}
e^{\beta t}\left\vert \Delta Y_{t}\right\vert ^{2}+\beta
\int_{t}^{T}e^{\beta s}\left\vert \Delta Y_{s-}\right\vert ^{2}ds 
&=&2\int_{t}^{T}e^{\beta s}\Delta Y_{s-}\Delta f(s)ds+2\int_{t}^{T}e^{\beta
s}\Delta Y_{s-}\Delta g(s)d\overleftarrow{B}_{s}\\&&+2\int_{t}^{T}e^{\beta
s}\Delta Y_{s-}d\left( \Delta K_{s}\right)  
-2\int_{t}^{T}\sum_{i=1}^{i=\infty }e^{\beta s}\Delta Y_{s-}\Delta
Z_{s}^{\left( i\right) }dH_{s}^{\left( i\right) }\\&&+\int_{t}^{T}e^{\beta
s}\left\vert \Delta g(s)\right\vert ^{2}ds-\int_{t}^{T}e^{\beta s}\Delta
Y_{s-}\sum_{i,j=1}^{\infty }\Delta Z_{s}^{\left( i\right) }\Delta
Z_{s}^{\left( j\right) }d\left[ H_{s}^{\left( i\right) },H_{s}^{\left(
j\right) }\right] .
\end{eqnarray*}%
Note that $\int_{0}^{t}e^{\beta s}\Delta Y_{s-}\Delta g(s)d\overleftarrow{B}%
_{s},$ $\int_{0}^{t}\sum_{i=1}^{i=\infty }e^{\beta s}\Delta Y_{s-}\Delta
Z_{s}^{\left( i\right) }dH_{s}^{\left( i\right) }$ $\forall $ $i\geq 1$ and%
\newline
$\int_{0}^{t}\sum_{i,j=1}^{\infty }e^{\beta s}\Delta Z_{s}^{\left( i\right)
}\Delta Z_{s}^{\left( j\right) }d\left( \left[ H_{s}^{\left( i\right)
},H_{s}^{\left( j\right) }\right] -\langle H_{s}^{\left( i\right)
},H_{s}^{\left( j\right) }\rangle \right) $ for $i\neq j$ are uniformly
integrable martingales. Since $\int_{t}^{T}e^{\beta s}\Delta Y_{s-}d\left(
\Delta K_{s}\right) \leq 0$, taking the mathematical expectation on bath
sides, we obtain%
\begin{eqnarray*}
&&\mathbb{E}e^{\beta t}\left\vert \Delta Y_{t}\right\vert ^{2}+\beta \mathbb{%
E}\int_{t}^{T}e^{\beta s}\left\vert \Delta Y_{s-}\right\vert ^{2}ds+\mathbb{E%
}\int_{t}^{T}\sum_{i=1}^{i=\infty }e^{\beta s}\left\vert \Delta
Z_{s}^{\left( i\right) }\right\vert ^{2}ds \\
&=&2\mathbb{E}\int_{t}^{T}e^{\beta s}\Delta Y_{s-}\Delta f(s)ds+\mathbb{E}%
\int_{t}^{T}e^{\beta s}\left\vert \Delta g(s)\right\vert ^{2}ds.
\end{eqnarray*}%
Firstly using $2ab\leq \epsilon _{1}a^{2}+\frac{b^{2}}{\epsilon _{1}}$ and
hypothesis $\left( H.1\right) ,$ we have%
\begin{equation*}
2\mathbb{E}\int_{t}^{T}e^{\beta s}\Delta Y_{s}\Delta f(s)ds\leq \epsilon _{1}%
\mathbb{E}\int_{t}^{T+K}e^{\beta s}\left\vert \Delta Y_{s-}\right\vert
^{2}ds+\left( \frac{c}{\epsilon _{1}}+cM\right) \mathbb{E}%
\int_{t}^{T+K}e^{\beta s}\left( \left\vert \Delta U_{s-}\right\vert
^{2}+\left\vert \left\vert \Delta V_{s}\right\vert \right\vert
_{l^{2}}^{2}\right) ds,
\end{equation*}%
and also%
\begin{equation*}
\mathbb{E}\int_{t}^{T}e^{\beta s}\left\vert \Delta g(s)\right\vert
^{2}ds\leq \left( c+cM\right) \mathbb{E}\int_{t}^{T+K}e^{\beta s}\left\vert
\Delta U_{s-}\right\vert ^{2}ds+\left( \alpha _{1}+\alpha _{2}M\right)
\mathbb{E}\int_{t}^{T+K}e^{\beta s}\left\vert \left\vert \Delta
V_{s}\right\vert \right\vert _{l^{2}}^{2}ds.
\end{equation*}
Then, we have
\begin{eqnarray*}
\mathbb{E}e^{\beta t}\left\vert \Delta Y_{t}\right\vert ^{2}&+&\left( \beta
-\epsilon _{1}\right) \mathbb{E}\int_{t}^{T}e^{\beta s}\left\vert \Delta
Y_{s-}\right\vert ^{2}ds+\mathbb{E}\int_{t}^{T+K}e^{\beta s}\left\vert
\left\vert \Delta Z_{s}\right\vert \right\vert _{l^{2}}^{2}ds \\
&\leq &\left( \frac{c}{\epsilon _{1}}+c+2cM\right) \mathbb{E}
\int_{t}^{T+K}e^{\beta s}\left\vert \Delta U_{s-}\right\vert ^{2}ds\\&+&\left(
\left( \alpha _{1}+\alpha _{2}M\right) +\left( \frac{c}{\epsilon _{1}}
+M\right) \right) \mathbb{E}\int_{t}^{T+K}e^{\beta s}\left\vert \left\vert
\Delta V_{s}\right\vert \right\vert _{l^{2}}^{2}ds,
\end{eqnarray*}%
we noting that $\Delta Y_{t}=\Delta Z_{t}=0,$ for all $t\in \left[ T,T+K%
\right] $%
\begin{eqnarray*}
\left( \beta -\epsilon _{1}\right) \mathbb{E}\int_{t}^{T+K}e^{\beta
s}\left\vert \Delta Y_{s-}\right\vert ^{2}ds&+&\mathbb{E}\int_{t}^{T+K}e^{%
\beta s}\left\vert \left\vert \Delta Z_{s}\right\vert \right\vert
_{l^{2}}^{2}ds 
\leq \left( \frac{c}{\epsilon _{1}}+c+2cM\right) \mathbb{E}%
\int_{t}^{T+K}e^{\beta s}\left\vert \Delta U_{s-}\right\vert ^{2}ds\\&+&\left(
\left( \alpha _{1}+\alpha _{2}M\right) +\left( \frac{c}{\epsilon _{1}}%
+M\right) \right) \mathbb{E}\int_{t}^{T+K}e^{\beta s}\left\vert \left\vert
\Delta V_{s}\right\vert \right\vert _{l^{2}}^{2}ds \\
&\leq &\left( \left( \alpha _{1}+\alpha _{2}M\right) +\left( \frac{c}{%
\epsilon _{1}}+cM\right) \right) [ \epsilon _{2}\mathbb{E}%
\int_{t}^{T+K}e^{\beta s}\left\vert \Delta U_{s-}\right\vert ^{2}ds\\&+&\mathbb{E%
}\int_{t}^{T+K}e^{\beta s}\left\vert \left\vert \Delta V_{s}\right\vert
\right\vert _{l^{2}}^{2}ds] ,
\end{eqnarray*}%
where $\epsilon _{2}=\frac{\frac{c}{\epsilon _{1}}+c+2cM}{\left( \alpha
_{1}+\alpha _{2}M\right) +\left( \frac{c}{\epsilon _{1}}+cM\right) }$.
Hence, if we choose $\epsilon _{0}=\epsilon _{1}$ satisfying $\hat{c}=\left(
\alpha _{1}+\alpha _{2}M\right) +\left( \frac{c}{\epsilon _{0}}+cM\right) <1$%
, chosse $\beta =\epsilon _{0}+\epsilon _{2}$. Then, we deduce%
\begin{equation*}
\mathbb{E}\int_{t}^{T+K}\epsilon _{2}e^{\beta s}\left\vert \Delta
Y_{s-}\right\vert ^{2}ds+\mathbb{E}\int_{t}^{T+K}e^{\beta s}\left\vert
\left\vert \Delta Z_{s}\right\vert \right\vert _{l^{2}}^{2}ds\leq \hat{c}%
\mathbb{E}\int_{t}^{T+K}e^{\beta s}\left( \epsilon _{2}\left\vert \Delta
U_{s-}\right\vert ^{2}+\left\vert \left\vert \Delta V_{s}\right\vert
\right\vert _{l^{2}}^{2}\right) ds.
\end{equation*}%
Thus, the mapping $\Phi $ is a strict contraction on $\mathcal{S}_{\mathcal{H%
}}^{2}\left( \left[ 0,T+K\right] ;\mathbb{R}^{d}\right) \times \mathcal{M}_{%
\mathcal{H}}^{2}\left( \left[ 0,T+K\right] ;l^{2}\right) $ and it has a
unique fixed point $\left( Y_{t},Z_{t}\right) \in \mathcal{S}_{\mathcal{H}%
}^{2}\left( \left[ 0,T+K\right] ;\mathbb{R}^{d}\right) \times \mathcal{M}_{%
\mathcal{H}}^{2}\left( \left[ 0,T+K\right] ;l^{2}\right) .$\eop

\subsection{Anticipated BDSDE.}

\qquad In this subsection we considere the anticipated BDSDE as follows

\begin{equation}
\left\{
\begin{array}{l}
Y_{t}=\xi _{T}+\int_{t}^{T}f(s,\Lambda _{s},\Lambda _{s}^{\phi ,\psi
})ds+\int_{t}^{T}g(s,\Lambda _{s},\Lambda _{s}^{\phi ,\psi })d\overleftarrow{%
B}_{s}-\sum_{i=1}^{\infty }\int_{t}^{T}Z_{s}^{\left( i\right)
}dH_{s}^{\left( i\right) },\text{ }t\in \left[ 0,T\right] , \\
\\
\left( Y_{t},Z_{t}\right) =\left( \eta _{t},\vartheta _{t}\right) ,\qquad
\qquad \qquad \qquad \qquad \qquad \qquad \qquad \qquad \qquad \qquad \qquad
t\in \left[ T,T+K\right] ,%
\end{array}%
\right.
\end{equation}%
where $f$ the generator, $\Lambda _{s}=\left( Y_{s-},Z_{s}\right) $ and $%
\Lambda _{s}^{\phi ,\psi }=\left( Y_{s+\phi \left( s\right) -},Z_{s+\psi
\left( s\right) }\right) .$

\begin{definition}
A solution of equation $(3.3)$ is a couple $(Y,Z)$ which belongs to the
space $S_{\mathcal{H}}^{2}([0,T+K],\mathbb{R})\times \mathcal{M}_{\mathcal{H}%
}^{2}([0,T+K];l^{2})$ and satisfies $(3.3)$.
\end{definition}

\subsubsection{\textbf{Existence and uniqueness of solutions.}}

\begin{theorem}
\label{th1} Assume that \textbf{(A), (B)} and \textbf{(H1)} are satisfied.
Then for given $\left( \eta _{t},\vartheta _{t}\right) \in \mathcal{B}_{%
\mathcal{H}}^{2}\left( \left[ T,T+K\right] ,%
\mathbb{R}
\right) $ l'equation $(3.3)$ has a unique solution $\left(
Y_{t},Z_{t}\right) \in \mathcal{B}_{\mathcal{H}}^{2}\left( \left[ 0,T+K%
\right] ,%
\mathbb{R}
\right) $.
\end{theorem}


\subsubsection{Proof of the existence and uniqueness result.}

\qquad Before we start proving equation $(3.3)$ has a unique solution with $%
f,$ $g$ independent on the value and the futur value of $\left( Y,Z\right) $%
, i.e., P-a.s.,$f\left( t,\omega ,y,z,\pi ,\zeta \right) =f\left( t,\omega
\right) $ and $g\left( t,\omega ,y,z,\pi ,\zeta \right) =g\left( t,\omega
\right) ,$ for any $\left( t,y,z,\pi ,\zeta \right) .$ More precisely, given
$f,$ $g$ such that%
\begin{eqnarray*}
E\left( \int_{0}^{T}\left\vert f\left( t\right) \right\vert ^{2}dt\right)
&<&\infty , \\
E\left( \int_{0}^{T}\left\vert g\left( t\right) \right\vert ^{2}dt\right)
&<&\infty .
\end{eqnarray*}%
Under the above assumptation on $f,$ $g,$ $\xi.$

\begin{proposition}
\label{Proposition1} Given $\xi\in \mathbb{L}^{2}\left( \mathcal{H}%
_{T}\right) ,$ there exists a unique couple of processes $\left(
Y_{t},Z_{t}\right) \in \mathcal{B}_{\mathcal{H}}^{2}\left( \left[ 0,T+K%
\right] ,%
\mathbb{R}
\right) ,$ to solve the following BDSDEs$,$%
\begin{equation*}
\left\{
\begin{array}{l}
Y_{t}=\xi _{T}+\int_{t}^{T}f(s)ds+\int_{t}^{T}g(s)d\overleftarrow{B}%
_{s}-\sum_{i=1}^{\infty }\int_{t}^{T}Z_{s}^{\left( i\right) }dH_{s}^{\left(
i\right) },\quad \ \ t\in \left[ 0,T\right] , \\
\left( Y_{t},Z_{t}\right) =\left( \eta _{t},\vartheta _{t}\right) ,\qquad
\qquad \qquad \qquad \qquad \qquad \qquad \qquad \ t\in \left[ T,T+K\right] .%
\end{array}%
\right.
\end{equation*}%
%
%
%
%
%
%
%
%
%
%
%
%
%
%
%
%
%
%
%
%
%
%
%
%
%
%
%
%
%
%
%
%
%
%
%
%
%
%
%
%
%
%
%
%
%
%
%
%
%
%
%
%
%
%
%
%
%
%
%
%
\end{proposition}

\noindent \textbf{Proof:} We consider the following filtration

\begin{equation*}
\mathcal{G}_{t}:=\mathcal{F}_{t}^{L}\vee \mathcal{F}_{T+K}^{B},
\end{equation*}%
and the $\mathcal{G}_{t}$ square integrable martingale%
\begin{equation*}
M_{t}=\mathbb{E}^{\mathcal{G}_{t}}\left( \xi
_{T}+\int_{t}^{T}f(s)ds+\int_{t}^{T}g(s)d\overleftarrow{B}_{s}\right) ,\text{%
\quad }t\in \left[ 0,T\right] .
\end{equation*}%
Thank's to the prédictable representation property in Nualart et al. [5] yields that there
existe $Z\in \mathcal{M}_{\mathcal{G}}^{2}\left( \left[ 0,T\right]
;l^{2}\right) $ such that%
\begin{equation*}
M_{t}=M_{0}+\sum_{i=1}^{i=\infty }\int_{0}^{t}Z_{s}^{\left( i\right)
}dH_{s}^{\left( i\right) },
\end{equation*}%
hence%
\begin{equation*}
M_{T}=M_{t}+\sum_{i=1}^{i=\infty }\int_{t}^{T}Z_{s}^{\left( i\right)
}dH_{s}^{\left( i\right) }.
\end{equation*}%
Let%
\begin{eqnarray*}
Y_{t} &=&M_{t}-\int_{t}^{T}f(s)ds-\int_{t}^{T}g(s)d\overleftarrow{B}_{s}, \\
&=&\mathbb{E}^{\mathcal{G}_{t}}\left( \xi
_{T}+\int_{t}^{T}f(s)ds+\int_{t}^{T}g(s)d\overleftarrow{B}_{s}\right) , \\
&=&M_{T}-\sum_{i=1}^{i=\infty }\int_{t}^{T}Z_{s}^{\left( i\right)
}dH_{s}^{\left( i\right) }-\int_{0}^{t}f(s)ds-\int_{0}^{t}g(s)d%
\overleftarrow{B}_{s},
\end{eqnarray*}%
from which, we deduce that%
\begin{equation*}
Y_{t}=\xi+\int_{t}^{T}f(s)ds+\int_{t}^{T}g(s)d\overleftarrow{B}%
_{s}-\sum_{i=1}^{i=\infty }\int_{t}^{T}Z_{s}^{\left( i\right)
}dH_{s}^{\left( i\right) }.
\end{equation*}%
Now by the same procedure of Xu [16], we can obtain the uniqueness and $%
\mathcal{H}_{t}-$ measurable of $Y_{t}$ and $Z_{t}$.\eop

We are now in a position to give the proof of Theorem 3.2.

\noindent \textbf{Proof:} Let $\mathcal{S}_{\mathcal{H}}^{2}\left( \left[
0,T+K\right] ;\mathbb{R}^{d}\right) \times \mathcal{M}_{\mathcal{H}%
}^{2}\left( \left[ 0,T+K\right] ;l^{2}\right) $ endowed with the norm%
\begin{equation*}
\left\vert \left\vert \left( Y,Z\right) \right\vert \right\vert _{\beta
}=\left( \mathbb{E}\left[ \int_{0}^{T+K}e^{\beta s}\left( \left\vert
Y_{s-}\right\vert ^{2}ds+\sum_{i=1}^{i=\infty }\left\vert Z_{s}^{\left(
i\right) }\right\vert ^{2}\right) ds\right] \right) ^{\frac{1}{2}}.
\end{equation*}%
Let we consider the following mapping:%
\begin{equation*}
\Phi :\mathcal{S}_{\mathcal{H}}^{2}\left( \left[ 0,T+K\right] ;\mathbb{R}%
\right) \times \mathcal{M}_{\mathcal{H}}^{2}\left( \left[ 0,T+K\right]
;l^{2}\right) \rightarrow \mathcal{S}_{\mathcal{H}}^{2}\left( \left[ 0,T+K%
\right] ;\mathbb{R}\right) \times \mathcal{M}_{\mathcal{H}}^{2}\left( \left[
0,T+K\right] ;l^{2}\right) ,
\end{equation*}%
where the couple $\left( Y_{t},Z_{t}\right) _{T\leq t\leq T+K}\in \mathcal{S}%
_{\mathcal{H}}^{2}\left( \left[ 0,T+K\right] ;\mathbb{R}\right) \times
\mathcal{M}_{\mathcal{H}}^{2}\left( \left[ 0,T+K\right] ;l^{2}\right) $ such
that

$\left( Y_{t},Z_{t}\right) _{T\leq t\leq T+K}=\left( \eta _{t},\vartheta
_{t}\right) $ and satisfies the equation $(E^{\xi ,f^{\phi ,\psi },g^{\phi
,\psi }})$. Thanks to Proposition $\left( 3.1\right) $, the mapping $\Phi $
is well defined. Let $\left( Y_{t},Z_{t}\right) $ and $\left( \tilde{Y}_{t},%
\tilde{Z}_{t}\right) $ be two solution of $(3.3)$ such that $\left(
Y_{t},Z_{t}\right) =\Phi \left( y_{t-},z_{t}\right) $ and $\left( \tilde{Y}%
_{t},\tilde{Z}_{t}\right) =\Phi \left( \tilde{y}_{t-},\tilde{z}_{t}\right) $.

For $\beta \in
\mathbb{R}
$. The couple $\left( \Delta Y_{t},\Delta Z_{t}\right) $ solve the ABDSDEs
with teugles martingale%
\begin{equation*}
\left\{
\begin{array}{l}
\Delta Y_{t}=\int_{t}^{T}\Delta f(s)ds+\int_{t}^{T}\Delta g(s)d%
\overleftarrow{B}_{s}-\sum_{i=1}^{i=\infty }\int_{t}^{T}\Delta Z_{s}^{\left(
i\right) }dH_{s}^{\left( i\right) },\text{ }t\in \left[ 0,T\right] , \\
\\
\left( \Delta Y_{t},\Delta Z_{t}\right) =\left( 0,0\right) ,\qquad \qquad
\qquad \qquad \qquad \qquad \qquad \qquad t\in \left[ T,T+K\right] .%
\end{array}%
\right.
\end{equation*}%
where for a function $h\in \left\{ f,g\right\} $, $\Delta h(s)=h(s,\theta
_{s},\theta _{s}^{\phi ,\psi })-h(s,\tilde{\theta}_{s},\tilde{\theta}%
_{s}^{\phi ,\psi }),$ $\theta _{s}=\left( y_{s-},z_{s}\right) ,$ $\theta
_{s}^{\phi ,\psi }=\left( y_{s+\phi \left( s\right) -},z_{s+\psi \left(
s\right) }\right) ,$ $\tilde{\theta}_{s}=\left( \tilde{y}_{s-},\tilde{z}%
_{s}\right) $, $\tilde{\theta}_{s}^{\phi ,\psi }=\left( \tilde{y}_{s+\phi
\left( s\right) -},\tilde{z}_{s+\psi \left( s\right) }\right) $ and $\Delta
\Psi _{s}=\Psi _{s}-\tilde{\Psi}_{s}$. Applying Itô's formula to $e^{\beta t}\left\vert \Delta Y_{t}\right\vert ^{2}$, we obtain%
\begin{eqnarray*}
e^{\beta t}\left\vert \Delta Y_{t}\right\vert ^{2}+\beta
\int_{t}^{T}e^{\beta s}\left\vert \Delta Y_{s-}\right\vert ^{2}ds 
&=&2\int_{t}^{T}e^{\beta s}\Delta Y_{s-}\Delta f(s)ds+2\int_{t}^{T}e^{\beta
s}\Delta Y_{s-}\Delta g(s)d\overleftarrow{B}_{s} \\
&&-2\int_{t}^{T}\sum_{i=1}^{i=\infty }e^{\beta s}\Delta Y_{s-}\Delta
Z_{s}^{\left( i\right) }dH_{s}^{\left( i\right) }+\int_{t}^{T}e^{\beta
s}\left\vert \Delta g(s)\right\vert ^{2}ds \\
&&-\int_{t}^{T}\sum_{i,j=1}^{\infty }e^{\beta s}\Delta Z_{s}^{\left(
i\right) }\Delta Z_{s}^{\left( j\right) }d\left[ H_{s}^{\left( i\right)
},H_{s}^{\left( j\right) }\right] ,
\end{eqnarray*}%
note that $\int_{0}^{t}e^{\beta s}\Delta Y_{s-}\Delta g(s)d\overleftarrow{B}%
_{s},$ $\int_{0}^{t}\sum_{i=1}^{i=\infty }e^{\beta s}\Delta Y_{s-}\Delta
Z_{s}^{\left( i\right) }dH_{s}^{\left( i\right) }$ $\forall $ $i\geq 1$ and%
\newline
$\int_{0}^{t}\sum_{i,j=1}^{\infty }e^{\beta s}\Delta Z_{s}^{\left( i\right)
}\Delta Z_{s}^{\left( j\right) }d\left( \left[ H_{s}^{\left( i\right)
},H_{s}^{\left( j\right) }\right] -\langle H_{s}^{\left( i\right)
},H_{s}^{\left( j\right) }\rangle \right) $ for $i\neq j$ are uniformly
integrable martingales. Now taking the mathematical expectation on bath
sides, we obtain%
\begin{eqnarray*}
&&\mathbb{E}e^{\beta t}\left\vert \Delta Y_{t}\right\vert ^{2}+\beta \mathbb{%
E}\int_{t}^{T}e^{\beta s}\left\vert \Delta Y_{s-}\right\vert ^{2}ds+\mathbb{E%
}\int_{t}^{T}\sum_{i=1}^{i=\infty }e^{\beta s}\left\vert \Delta
Z_{s}^{\left( i\right) }\right\vert ^{2}ds \\
&=&2\mathbb{E}\int_{t}^{T}e^{\beta s}\Delta Y_{s-}\Delta f(s)ds+\mathbb{E}%
\int_{t}^{T}e^{\beta s}\left\vert \Delta g(s)\right\vert ^{2}ds.
\end{eqnarray*}%
Now by the same computation of Lipschitz coefficient for Anticipated
reflected BDSDEs in general case, we deduce that%
\begin{equation*}
\mathbb{E}\int_{t}^{T+K}e^{\beta s}\left( \epsilon _{2}\left\vert \Delta
Y_{s-}\right\vert ^{2}+\left\vert \left\vert \Delta Z_{s}\right\vert
\right\vert _{l^{2}}^{2}\right) ds\leq \varepsilon \mathbb{E}%
\int_{t}^{T+K}e^{\beta s}\left( \epsilon _{2}\left\vert \Delta
y_{s-}\right\vert ^{2}+\left\vert \left\vert \Delta z_{s}\right\vert
\right\vert _{l^{2}}^{2}\right) ds,
\end{equation*}%
where $0<\varepsilon <1.$ Thus, the mapping $\Phi $ is a strict contraction
on $\mathcal{S}_{\mathcal{H}}^{2}\left( \left[ 0,T+K\right] ;\mathbb{R}%
\right) \times \mathcal{M}_{\mathcal{H}}^{2}\left( \left[ 0,T+K\right]
;l^{2}\right) $ and it has a unique fixed point $\left( Y_{t},Z_{t}\right)
\in \mathcal{S}_{\mathcal{H}}^{2}\left( \left[ 0,T+K\right] ;\mathbb{R}%
\right) \times \mathcal{M}_{\mathcal{H}}^{2}\left( \left[ 0,T+K\right]
;l^{2}\right) .$

Finally we complete the proof of theorem 3.2.\eop


\end{document}